\documentclass{article}
\usepackage{amsmath,amssymb,amsthm,dsfont,hyperref,tikz}
\usepackage{geometry}
\allowdisplaybreaks
\begin{document}
    \title{On Correlation Functions of Vertex Operator Algebras Associated to Jordan Algebras}
    \author{Hongbo Zhao}
    \date{}
    \maketitle
    \newtheorem{lemma}{Lemma}
    \newtheorem{theorem}{Theorem}
    \newtheorem{proposition}{Proposition}
    \allowdisplaybreaks
    \section{Introduction}
    \par{
        Given a $\mathbb{Z}_{\geq 0}$ graded vertex operator algebra (VOA) $V=\bigoplus_{i=0}^{\infty}V_i$ with $dim(V_0)=1$, then $V_1$ has a structure of Lie algebra, with the operation given by $[a,b]=a(0)b, \;\forall a,b\in V_1$. The affine vertex algebra $V_r(\mathfrak{g})$ of level $r$ provides such an example, with $V_1=\mathfrak{g}$ as a Lie algebra.
    }
    \par{
         When $dim(V_0)=1$, $dim(V_1)=0$, then $V_2$ has a structure of commutative (but not necessarily associative) algebra, with the operation $a\circ b=a(1)b$.
        For moonshine VOA $V^{\natural}$ \cite{FLM} with $dim(V^{\natural}_0)=1,\;dim(V^{\natural}_1)=0$, the corresponding commutative algebra $V^{\natural}_2$ is the Griess algebra introduced by Griess in the study of the monster group \cite{Gr}. We call $V_2$ of a VOA $V$ with $dim(V_0)=1$, $dim(V_1)=0$ the Griess algebra of $V$. In \cite{Lam1} and \cite{Lam2}, Lam constructed vertex algebras whose Griess algebras are simple Jordan algebras of type $A,B,C$. In \cite{AM}, Ashihara and Miyamoto constructed for a Jordan algebra $\mathcal{J}$ of Type $B$, a family of VOAs $V_{\mathcal{J},r}$ parameterized by a complex number $r$, such that
        $(V_{\mathcal{J},r})_0=\mathbb{C}1,(V_{\mathcal{J},r})_1=\{0\}$, and
        $$
            (V_{\mathcal{J},r})_2\cong\mathcal{J}
        $$
        as Jordan algebra, and Lam's example for type $B$ Jordan algebra is a quotient of $V_{\mathcal{J},1}$.
        The main result of this paper is a formula about the genus zero correlation function of generating fields in $V_{\mathcal{J},r}$.
    }
    \par{A. Albert classified finite dimensional simple Jordan algebras over an algebraic closed field $\mathbb{F}$ with $char(\mathbb{F})=0$ \cite{A}. We have a brief review of simple Jordan algebra of type $B$, and we only consider the case $\mathbb{F}=\mathbb{C}$.
        Let $(\mathfrak{h},(\cdot,\cdot))$ be a finite dimensional vector space with a non-degenerate symmetric bilinear form $(\cdot,\cdot)$ and $dim(\mathfrak{h})=d$. Then $\mathfrak{h}\otimes \mathfrak{h}$ has an associative algebra  structure:
        $$
            (a\otimes b)(u\otimes v)=(b,u)a\otimes v,
        $$
        which induces a Jordan algebra structure on $\mathfrak{h}\otimes \mathfrak{h}$:
        $$
            x\circ y=\frac{1}{2}(xy+yx),\;\forall x,y\in \mathfrak{h}\otimes \mathfrak{h}.
        $$
        Let $\mathcal{J}$ be the Jordan subalgebra of $\mathfrak{h}\otimes \mathfrak{h}$ consists of symmetric tensors:
        $$
            \mathcal{J}\stackrel{def}{=}span\{L_{a,b}|a,b\in\mathfrak{h}\},\;L_{a,b}\stackrel{def}{=}a\otimes b+b\otimes a.
        $$
        Then $\mathcal{J}$ is essentially the Jordan algebra of symmetric $d\times d$ matrices over $\mathbb{C}$, which is called simple Jordan algebra of type $B$ according to Jacobson's notation\cite{Jac}.
    }

    \par{
        To state our result, we need to introduce some notations. Let $L_{a,b}(z)$ denote the vertex operator associated to $L_{a,b}$. Given $n$ vertex operators $L_{a_1,b_1}(z_1),\cdots ,L_{a_{n},b_n}(z_n)$, we have a corresponding sequence $T=(a_1,b_1)\cdots(a_n,b_n)$.
        Recall that a derangement of the $n$-element set $\{(a_1,b_1),\cdots,(a_n,b_n)\}$ is a permutation of $n$-elements without fix point \cite{FS}. Let $DR(T)$ denote the set of all derangements corresponding to $T$. For $\sigma\in DR(T)$ we decompose it as disjoint cycles£º
        $$
            \sigma=(C_1)\cdots(C_s),
        $$
        and we use
        $\sigma(i)$ to denote the image of each label $i$ under the action of $\sigma$. We also use the symbol
        $$
            c(\sigma)\stackrel{def}{=}s
        $$
        to denote the number of disjoint cycles in $\sigma$.
    }
    \par{
        Our main result is stated as follows:
		\begin{theorem}
			Given $n$ vertex operators $L_{a_1,b_1}(z_1),\cdots, L_{a_{n},b_n}(z_n)$, and $T=(a_1,b_1)\cdots(a_n,b_n)$. Then the corresponding genus zero correlation function is given by
			\begin{align}
				\langle1^{'},L_{a_1,b_1}(z_1)\cdots L_{a_n,b_n}(z_n)\cdot 1\rangle
=\sum_{\sigma\in DR(T)}\Gamma(\sigma,T)\Gamma(\sigma;Z)r^{c(\sigma)},\label{form1}
			\end{align}
			where the symbols are described as follows: Assume that $\sigma=(C_1)\cdots(C_s)=(k_{11}\cdots k_{1 t_1})\cdots (k_{s1}\cdots k_{s t_s})$, then
    \begin{align*}
        \Gamma(\sigma,T)\stackrel{def}{=}&2^{-s-n}\prod_{i=1}^{s}Tr(L_{a_{k_{i1}},b_{k_{i1}}}\cdots L_{a_{k_{it_i}},b_{k_{it_i}}}),\\
        \Gamma(\sigma;Z)\stackrel{def}{=}&\prod_{i=1}^n\frac{1}{(z_i-z_{\sigma(i)})^2},
    \end{align*}
    where $Tr(a\otimes b)$ is the trace of $a\otimes b \in \mathfrak{h}\otimes \mathfrak{h}$ given by:
    $$
        Tr(a\otimes b)=(a,b).
    $$
    \end{theorem}
	}
\par{
    We have a more general formula about the ``correlation function" of some two variable formal power series, which will be shown in Proposition 2, and we prove Theorem 1 as a corollary of this proposition.
}
\par{
    We give an example of Theorem 1 for the case $n=4$. The sequence is $T=(a_1,b_1)\cdots(a_4,b_4)$, and there are $9$ elements in $DR(T)$, given by:
    $$
        (12)(34),(13)(24),(14)(23);(1234),(1243),(1324),(1342),(1423),(1432).
    $$
    Then Theorem 1 gives
    \begin{align*}
        &\langle 1^{'},L_{a_1,b_1}(z_1)L_{a_2,b_2}(z_2)L_{a_3,b_3}(z_3)L_{a_4,b_4}(z_4)1\rangle\\
        =&\frac{Tr(L_{a_1,b_1}L_{a_2,b_2})Tr(L_{a_3,b_3}L_{a_4,b_4})r^2}{64(z_1-z_2)^4(z_3-z_4)^4}\notag\\&+\frac{Tr(L_{a_1,b_1}L_{a_3,b_3})Tr(L_{a_2,b_2}L_{a_4,b_4})r^2}{64(z_1-z_3)^4(z_2-z_4)^4}+\frac{Tr(L_{a_1,b_1}L_{a_4,b_4})Tr(L_{a_2,b_2}L_{a_3,b_3})r^2}{64(z_1-z_4)^4(z_2-z_3)^4}\\
        &+\frac{Tr(L_{a_1,b_1}L_{a_2,b_2}L_{a_3,b_3}L_{a_4,b_4})r}{32(z_1-z_2)^2(z_2-z_3)^2(z_3-z_4)^2(z_4-z_1)^2}+\frac{Tr(L_{a_1,b_1}L_{a_2,b_2}L_{a_4,b_4}L_{a_3,b_3})r}{32(z_1-z_2)^2(z_2-z_4)^2(z_4-z_3)^2(z_3-z_1)^2}\\
        &+\frac{Tr(L_{a_1,b_1}L_{a_3,b_3}L_{a_2,b_2}L_{a_4,b_4})r}{32(z_1-z_3)^2(z_3-z_2)^2(z_2-z_4)^2(z_4-z_1)^2}+\frac{Tr(L_{a_1,b_1}L_{a_3,b_3}L_{a_4,b_4}L_{a_2,b_2})r}{32(z_1-z_3)^2(z_3-z_4)^2(z_4-z_2)^2(z_2-z_1)^2}\\
        &+\frac{Tr(L_{a_1,b_1}L_{a_4,b_4}L_{a_2,b_2}L_{a_3,b_3})r}{32(z_1-z_4)^2(z_4-z_2)^2(z_2-z_3)^2(z_3-z_1)^2}+\frac{Tr(L_{a_1,b_1}L_{a_4,b_4}L_{a_3,b_3}L_{a_2,b_2})r}{32(z_1-z_4)^2(z_4-z_3)^2(z_3-z_2)^2(z_2-z_1)^2}.
    \end{align*}
}
\par{
    When $r=1$, the correlation function in Theorem 1 is the same as correlation function of vertex operators $\frac{1}{2}:a_1(z_1)b_1(z_1):,\cdots,\frac{1}{2}:a_n(z_n)b_n(z_n):$ in Heisenberg VOA $S(\hat{\mathfrak{h}}_{-})$. The theorem can be proved using Wick's theorem (see Theorem 3.3 in \cite{Kac}).
}
\par{
    When $dim(\mathcal{J})=1$, the VOA $V_{\mathcal{J},r}$ is the same as the VOA associated to Virasoro algebra with central charge equal to $r$, and the conformal vector is $\omega=L_{e,e}$, $(e,e)=1$. In this case our formula agrees with the result in \cite{Hur} (Theorem 2.3), where the genus zero correlation function for Virasoro field is computed:
    \begin{align*}
				\langle 1^{'},\omega(z_1)\cdots \omega(z_n)\cdot 1\rangle=\sum_{\sigma\in DR(T)}\Gamma(\sigma;Z)(\frac{r}{2})^{c(\sigma)}=\sum_{\sigma\in DR(T)}(\prod_{i=1}^n\frac{1}{(z_i-z_{\sigma(i)})^2})(\frac{r}{2})^{c(\sigma)}.
	\end{align*}
    When $n=4$, the result is:
    \begin{align*}
        &\langle 1^{'},\omega(z_1)\omega(z_2)\omega(z_3)\omega(z_4)1\rangle\\
        =&\frac{r^2}{4(z_1-z_2)^4(z_3-z_4)^4}+\frac{r^2}{4(z_1-z_3)^4(z_2-z_4)^4}+\frac{r^2}{4(z_1-z_4)^4(z_2-z_3)^4}\\
        &+\frac{r}{2(z_1-z_2)^2(z_2-z_3)^2(z_3-z_4)^2(z_4-z_1)^2}+\frac{r}{2(z_1-z_2)^2(z_2-z_4)^2(z_4-z_3)^2(z_3-z_1)^2}\\
        &+\frac{r}{2(z_1-z_3)^2(z_3-z_2)^2(z_2-z_4)^2(z_4-z_1)^2}+\frac{r}{2(z_1-z_3)^2(z_3-z_4)^2(z_4-z_2)^2(z_2-z_1)^2}\\
        &+\frac{r}{2(z_1-z_4)^2(z_4-z_2)^2(z_2-z_3)^2(z_3-z_1)^2}+\frac{r}{2(z_1-z_4)^2(z_4-z_3)^2(z_3-z_2)^2(z_2-z_1)^2}.
    \end{align*}
}
\par{
        Theorem 1 can also be viewed as an analogue of Theorem 2.3.1 in \cite{FZ}, where the genus zero correlation function of the generating fields in $V_r(\mathfrak{g})$ is computed:
}
    \par{
        \begin{theorem}
			The genus zero correlation function of fields $a_1(z_1),\cdots,a_n(z_n)$, $a_i\in\mathfrak{g}$ is given by the formula:
            $$
                \langle1^{'},a_1(z_1)\cdots a_n(z_n)1\rangle=\sum_{\sigma=(C_1)\cdots(C_s)\in DR(T)}\prod_{i=1}^s f_{C_i}(a_1,\cdots,a_n;z_1,\cdots,z_n)(-r)^s,
            $$
            where
            $$
                f_{C_i}(a_1,\cdots,a_n;z_1,\cdots,z_n)\stackrel{def}{=}\frac{Tr(a_{k_1}\cdots a_{k_t})}{(z_{k_1}-z_{k_2})\cdots(z_{k_{t-1}}-z_{k_t})(z_{k_{t}}-z_{k_1})}.
            $$
            The trace here is normalized so that $\forall x,y\in\mathfrak{g}$, $Tr(xy)$ equals to the normalized Killing form of $\mathfrak{g}$.
        \end{theorem}
    }
    \par{
        The content of this paper is organized as follows: In section 2 we have an overview of the construction in \cite {AM} and some results in \cite{NS}, then we introduce some two variable formal power series whose commutation relations are computed. In section 3, we introduce diagrams and some other related notations, and section 4 is devoted to the proof of Theorem 1.
    }
    \par{
        \textbf{Acknowledgement.} I would like to express my gratitude to my advisor, prof. Y. Zhu, for suggesting this problem.
    }
    \section{Vertex Operator Algebras Associated to Jordan Algebras of Type $B$}
    \par{
        In this section, first we give a brief overview of the vertex operator algebra $V_{\mathcal{J},r}$ constructed in \cite{AM}. Then, we introduce two-variable formal power series $L_{a,b}(z,w),L_{a,b}^{\pm\pm}(z,w)$ and compute some commutation relations.
    }
    \par{
        Let $\mathfrak{h}$ be a finite dimensional vector space with a non-degenerate symmetric bilinear form $(\cdot,\cdot)$. Assume $dim(\mathfrak{h})=d$, then the Heisenberg Lie algebra associated to $(\mathfrak{h},(\cdot,\cdot))$ is:
        \begin{align*}
            \hat{\mathfrak{h}}=L(\mathfrak{h})\oplus \mathbb{C}c,\;\;
            L(\mathfrak{h})=\mathfrak{h}\otimes\mathbb{C}[t,t^{-1}],
        \end{align*}
        with the Lie bracket  given by:
        \begin{align}
            [a(m),b(n)]=m(a,b)\delta_{m+n}c,\;\;
            [x,c]=0,\;\;\forall x\in \hat{\mathfrak{h}},\notag
        \end{align}
        where $a(m)=at^m\in L(\mathfrak{h})$. Noting that
        $$
            \hat{\mathfrak{h}}_{-}=\mathfrak{h}\otimes\mathbb{C}t^{-1}[t^{-1}]
        $$
        is a commutative Lie subalgebra, and
        the Fock space $S(\hat{\mathfrak{h}}_{-})=U(\hat{\mathfrak{h}}_{-})$  has a vertex operator algebra structure \cite{FLM}.

    \par{Take the following subspace $\mathfrak{L}$ of the universal enveloping algebra $U(\hat{\mathfrak{h}})$ spanned by quadratic elements:
        $$
            \mathfrak{L}\stackrel{def}{=}span\{a(m)b(n)|a,b\in \mathfrak{h},m,n\in\mathbb{Z}\}.
        $$
        Since $[a(m),b(n)]=m(a,b)\delta_{m+n}c$ so $c\in \mathfrak{L}$.
    }
    \par{
        Define a new Lie bracket '$[x,y]_{new}$' on $\mathfrak{L}$ by
        \begin{align}
            [x,y]_{new}\stackrel{def}{=}\frac{1}{c}[x,y],\;\forall x,y\in\mathfrak{L}.\label{form3}
        \end{align}
        The identity (\ref{form4}) below shows that the right hand side is still in $\mathfrak{L}$.
    }

\par{
    For convenience we use the notation:
    $$
        L_{a,b}(m,n)\stackrel{def}{=}\frac{1}{2}:a(m)b(n):,
    $$
    where
    \begin{align*}
        	:a(m)b(n):=
        	\begin{cases}
        		b(n)a(m),m\geq n,\\
        		a(m)b(n),else.
        	\end{cases}
        \end{align*}
        Let
        $$
            \mathfrak{B}\stackrel{def}{=}span\{L_{a,b}(m,n)|a,b\in \mathfrak{h},m,n\in\mathbb{Z}\}.
        $$
        Noting that we have a split of $\mathfrak{L}$:
        \begin{align*}
        	\mathfrak{L}=\mathfrak{B}\oplus\mathbb{C}c.
        \end{align*}
}
       \par{
        For two elements $L_{a,b}(m,n),L_{u,v}(p,q)\in \mathfrak{B}$, we have
        \begin{align}
         	&[L_{a,b}(m,n),L_{u,v}(p,q)]_{new}\notag\\
         	=&\frac{1}{4c}[a(m)b(n),u(p)v(q)]\notag\\
         	=&\frac{1}{4c}[a(m) b(n),u(p)] v(q)+\frac{1}{4c}u(p) [a(m)b(n),v(q)]\notag\\
         	=&\frac{1}{4c}a(m)[b(n),u(p)] v(q)+\frac{1}{4c}[a(m),u(p)] b(n)v(q)+\frac{1}{4c}u(p) a(m)[b(n),v(q)]+\frac{1}{4c}u(p)[a(m),v(q)] b(n)\notag\\
         	=&\frac{1}{4}n\delta_{n+p}(b,u)a(m) v(q)+\frac{1}{4}m\delta_{m+p}(a,u) b(n) v(q)\notag\\&+\frac{1}{4}n\delta_{n+q}(b,v)u(p) a(m)+\frac{1}{4}m\delta_{m+q}(a,v)u(p) b(n)\in \mathfrak{L}. \label{form4}
        \end{align}
    }
    \par{
        Decompose $\mathfrak{B}$ into:
        $$
            \mathfrak{B}=\mathfrak{B}_{+}\oplus\mathfrak{B}_{-},
        $$
        where
        \begin{align}
            &\mathfrak{B}_{-}\stackrel{def}{=}span\{L_{a,b}(m,n)|m,n<0\},\;\;\;
            \mathfrak{B}_{+}\stackrel{def}{=}span\{L_{a,b}(m,n)|n\geq 0\;\text{or}\; m\geq 0\}.\notag
        \end{align}
        Then we have a decomposition of $\mathfrak{L}$:
        \begin{align}
            &\mathfrak{L}_{-}=\mathfrak{B}_{-},\notag,\;\;\;\mathfrak{L}_{+}=\mathfrak{B}_{+}\oplus\mathbb{C}c.\notag
        \end{align}
    }
    \par{
       	Define a 1-dimensional $\mathfrak{L}_{+}$ module $\mathbb{C}1$ :
        \begin{align}
        	&x 1=0,\;\;\forall x\in\mathfrak{B}_{+},\notag\;\;\;c 1=r1.\notag
        \end{align}
        Then by induction from $U(\mathfrak{L}_{+})$ to $U(\mathfrak{L})$, we have a $U(\mathfrak{L})$ module $M_r$:
        \begin{align*}
        	M_r\stackrel{def}{=}&U(\mathfrak{L})\otimes_{U(\mathfrak{L}_{+})}\mathbb{C}1\\
        \cong&U(\mathfrak{L}_{-})1\\
        =&span\{L_{a_1,b_1}(-m_1,-n_1)\cdots L_{a_k,b_k}(-m_k,-n_k)\cdot 1| m_i,n_i\in \mathbb{Z}_{\geq 1},a_i,b_i\in \mathfrak{h}\}.
        \end{align*}
    }
    \par{
		For $a,b\in\mathfrak{h}$ define the following operator $L_{a,b}(l)$:
		\begin{align*}
        L_{a,b}(l)\stackrel{def}{=}&\sum_{k\in\mathbb{Z}}L_{a,b}(-k+l-1,k),
        \end{align*}
        and field $L_{a,b}(z)$ by:
        \begin{align*}
        L_{a,b}(z)\stackrel{def}{=}&\sum_{l\in\mathbb{Z}}L_{a,b}(l)z^{-l-1}.
		\end{align*}
        It is proved in \cite{AM} that these fields are mutually local.
}
\par{
		So these mutually local fields generate a vertex algebra, which is denoted by $V_{\mathcal{J},r}$:
		$$
			V_{\mathcal{J},r}{=}span\{L_{a_1,b_1}(m_1)\cdots L_{a_k,b_k}(m_k)\cdot 1| m_i\in \mathbb{Z},a_i,b_i\in \mathfrak{h}\}.
		$$
        It is proved in \cite{NS} that $V_{\mathcal{J},r}=M_r$ holds, the Virasoro element of $V_{\mathcal{J},r}$ is
        $$
            \omega=\sum_k L_{e_k,e_k}(-1,-1)\cdot 1,
        $$
        where $\{e_k\}$ is an orthonormal basis of $\mathfrak{h}$, and $\omega(1)$ gives a graduation of $V_{\mathcal{J},r}$:
        $$
            V_{\mathcal{J},r}=\sum_{k\geq 0}(V_{\mathcal{J},r})_{k}.
        $$
        The Griess algebra $(V_{\mathcal{J},r})_2$ is isomorphic to $\mathcal{J}$:
        $$
            L_{a,b}(-1,-1)\cdot 1\mapsto L_{a,b}\stackrel{def}{=}a\otimes b+b\otimes a.
        $$
}
\par{
    We introduce two-variable formal power series $L_{a,b}^{--}(z,w),\;L_{a,b}^{-+}(z,w),\;L_{a,b}^{+-}(z,w),\;L_{a,b}^{++}(z,w)$, and $L_{a,b}(z,w)$ whose coefficients are in $\mathfrak{L}$:
    \begin{align*}
    &L_{a,b}(z,w)\stackrel{def}{=}\sum_{m,n\in\mathbb{Z}} L_{a,b}(m,n)z^{-m-1}w^{-n-1},\\
    &L_{a,b}^{--}(z,w)\stackrel{def}{=}\sum_{m,n<0} L_{a,b}(m,n)z^{-m-1}w^{-n-1},\;
    L_{a,b}^{-+}(z,w)\stackrel{def}{=}\sum_{m<0,n\geq 0} L_{a,b}(m,n)z^{-m-1}w^{-n-1},\\
    &L_{a,b}^{+-}(z,w)\stackrel{def}{=}\sum_{m\geq 0,n<0} L_{a,b}(m,n)z^{-m-1}w^{-n-1},\;
    L_{a,b}^{++}(z,w)\stackrel{def}{=}\sum_{m,n\geq 0} L_{a,b}(m,n)z^{-m-1}w^{-n-1}.
    \end{align*}
    It is easy to see that
    \begin{align}
        &L_{a,b}^{--}(z,w)=L_{b,a}^{--}(w,z),\;L_{a,b}^{++}(z,w)=L_{b,a}^{++}(w,z),\;L_{a,b}^{+-}(z,w)=L_{b,a}^{-+}(w,z),\;\notag\\
        &L_{a,b}(z,w)=L_{a,b}^{--}(z,w)+L_{a,b}^{-+}(z,w)+L_{a,b}^{+-}(z,w)+L_{a,b}^{++}(z,w).\label{form17}
    \end{align}
    and
    $$
        L_{a,b}(z)=L_{a,b}(z,z)
    $$
    as fields on $V_{\mathcal{J},r}$.
}
\par{
    Then we have the following closed formula for the commutation relations that we need later:
    \begin{proposition}
    As formal power series with coefficients taking value in $End(V_{\mathcal{J},r})$, we have
    \begin{align}
        [L_{a,b}^{++}(x,y),L_{u,v}^{--}(z,w)]=&\frac{1}{2}(b,v)\iota_{y,w}(y-w)^{-2}L_{u,a}^{-+}(z,x)\notag+\frac{1}{2}(a,v)\iota_{x,w}(x-w)^{-2}L_{u,b}^{-+}(z,y)&\notag\\&+\frac{1}{2}(a,u)\iota_{x,z}(x-z)^{-2}
        L_{v,b}^{-+}(w,y)+\frac{1}{2}(b,u)\iota_{y,z}(y-z)^{-2}L_{v,a}^{-+}(w,x)\notag\\&+
        \frac{1}{4}r(a,u)(b,v)\iota_{x,z}(x-z)^{-2}\iota_{y,w}(y-w)^{-2}\notag\\&+\frac{1}{4}r(a,v)(b,u)\iota_{x,w}(x-w)^{-2}\iota_{y,z}(y-z)^{-2},\label{form5}\\
        [L_{a,b}^{++}(x,y),L_{u,v}^{-+}(z,w)]=&\frac{1}{2}(a,u)L_{b,v}^{++}(y,w)\iota_{x,z}(x-z)^{-2}+\frac{1}{2}(b,u)L_{a,v}^{++}(x,w)\iota_{y,z}(y-z)^{-2},\label{form6}\\
        [L_{a,b}^{++}(x,y),L_{u,v}^{++}(z,w)]=&0.\notag
    \end{align}
    \end{proposition}
    \textbf{Proof.} We only prove (\ref{form5}), and others are obtained in a similar way.
        Consider the following formal power series with coefficients in $U(\hat{\mathfrak{h}})$
        \begin{align*}
            a_{-}(z)\stackrel{def}{=}&\sum_{k<0} a(m)z^{-k-1},\;\;\;
            a_{+}(z)\stackrel{def}{=}\sum_{k\geq 0} a(m)z^{-k-1},\\
            a(z)\stackrel{def}{=}&\sum_k a(m)z^{-k-1}=a_{-}(z)+a_{+}(z),\;\;\forall a\in\mathfrak{h}.
        \end{align*}
        A direct computation shows that
        \begin{align}
            [a_{+}(z),b_{-}(w)]=(a,b)c\iota_{z,w}\frac{1}{(z-w)^2}\label{form2}.
        \end{align}
    }
    \par{
        We view $a_{\pm}(x)b_{\pm}(y)$ as formal power series with coefficients in $\mathfrak{L}$, then we have
        \begin{align*}
            &L^{++}_{a,b}(z,w)=\frac{1}{2}a_{+}(z)b_{+}(w),\;L^{-+}_{a,b}(z,w)=\frac{1}{2}a_{-}(z)b_{+}(w),\\&L^{+-}_{a,b}(z,w)=\frac{1}{2}b_{-}(w)a_{+}(z),\;L^{--}_{a,b}(z,w)=\frac{1}{2}a_{-}(z)b_{-}(w).
        \end{align*}
        So \ref{form5} is proved by the following computation:
    \begin{align}
        &[a_{+}(x)b_{+}(y),u_{-}(z)v_{-}(w)]_{new}\notag\\=&\frac{1}{c}[a_{+}(x)b_{+}(y),u_{-}(z)v_{-}(w)]\notag\\
        =&(b,v)\iota_{y,w}(y-w)^{-2}u_{-}(z)a_{+}(x)+(a,v)\iota_{x,w}(x-w)^{-2}u_{-}(z)b_{+}(y)\notag\\
        &+(a,u)\iota_{x,z}(x-z)^{-2}b_{+}(y)v_{-}(w)+(b,u)\iota_{y,z}(y-z)^{-2}a_{+}(x)v_{-}(w)\notag\\
        =&(b,v)\iota_{y,w}(y-w)^{-2}u_{-}(z)a_{+}(x)+(a,v)\iota_{x,w}(x-w)^{-2}u_{-}(z)b_{+}(y)\notag\\
        &+(a,u)\iota_{x,z}(x-z)^{-2}v_{-}(w)b_{+}(y)+(b,u)\iota_{y,z}(y-z)^{-2}v_{-}(w)a_{+}(x)\notag\\
        &+(a,u)(b,v)c\iota_{x,z}(x-z)^{-2}\iota_{y,w}(y-w)^{-2}\notag\\&+(a,v)(b,u)c\iota_{x,w}(x-w)^{-2}\iota_{y,z}(y-z)^{-2}.\notag
    \end{align}

}
\section{Diagrams, Derangements, and Some Necessary Notations}
    \par{
        In this section we introduce a notation called diagram over a sequence $T=(a_1,b_1)\cdots(a_n,b_n)$, denoted by $D(T)$, and it will be shown that there is a surjective map from $D(T)$ to $DR(T)$. We also introduce the notation of diagram over $T$ compatible with sign $\epsilon$, denoted by $D(T^{\epsilon})$, and the two obvious operations of diagram turn out to have correspondences with the commutation relations, which will be used to prove Theorem 1.
    }
    \par{
        Let $1^{'}$ be the unique element in $V_{\mathcal{J},r}^{*}$ satisfying $\langle1^{'},1\rangle$=1, $\langle1^{'},v\rangle=0,\forall v\in (V_{\mathcal{J},r})_{i\geq 1}$. Then the genus zero correlation function of $L_{a_1,b_1}(z_1),\cdots, L_{a_{n},b_n}(z_n)$ given by
        $$
            \langle 1^{'}, L_{a_1,b_1}(z_1)\cdots L_{a_{n},b_n}(z_n)1\rangle,
        $$
        which is a power series in $\mathbb{C}[[z_1,z^{-1}_1\cdots z_n,z^{-1}_n]]$. According to \cite{FLM}, this formal power series converges to a rational function in
        the domain $|z_1|>\cdots>|z_n|$.
    }
    \par{
        Given the sequence
        $T=(a_1,b_1)\cdots(a_n,b_n)$, a sign on $T$ is a function $\epsilon:\{a_1,b_1,\cdots,a_n,b_n\}\rightarrow \{+,-\}$. If $\epsilon(a_i)=\epsilon_i,\epsilon(b_i)=\delta_i$, we also write
        $$
            T^{\epsilon}=(a_1^{\epsilon_1},b_1^{\delta_1})\cdots(a_n^{\epsilon_n},b_n^{\delta_n}),
        $$
        and we call $(a_i^{\epsilon_i},b_i^{\delta_i})$ the $i$-th  pair of $T^{\epsilon}$.
    }
    \par{
        A diagram over the sequence $T=(a_1,b_1)\cdots(a_n,b_n)$ is a graph, with vertex set $V=\{a_1,b_1,\cdots,a_n,b_n\}$, and edge set $E$ consists of unordered pairs $\{u,v\},\;u,v\in V$ satisfying:
        \begin{description}
            \item[(1).] $|E|=n$.
            \item[(2).] $\{a_i,b_i\},\{a_i,a_i\},\{b_i,b_i\}\notin E,\forall i=1,\cdots,n$.
            \item[(3).] Any two edges have no common point.
        \end{description}
        Denote the set of all diagrams over $T$ by $D(T)$.
        Observing that such a graph is always with $2n$ vertices and $n$ disjoint edges.
    }
    \par{
        A  diagram over the signed sequence $T^{\epsilon}=(a_1^{\epsilon_1},b_1^{\delta_1})\cdots(a_n^{\epsilon_n},b_n^{\delta_n})$ is a diagram over $T$ compatible with the sign $\epsilon$ in the following sense:
        \begin{description}
            \item[(4).] $\forall \{u,v\}\in E$, if $u=a_i$ or $b_i$, $v=a_j$ or $b_j$, $i<j$, then $\epsilon(u)=+,\epsilon(v)=-$.
        \end{description}
        Denote the set of all  diagrams over $T$ which are compatible with $\epsilon$ by $D(T^{\epsilon})$. It may happen that for some $\epsilon$, $D(T^{\epsilon})=\emptyset$. For example, if $T^{\epsilon}=(a_1^{+},b_1^{+})(a_2^{-},b_2^{-})(a_3^{-},b_3^{-})(a_4^{+},b_4^{+})$, then $D(T^{
        \epsilon})=\emptyset$.
    }
    \par{

        Given a diagram in $\coprod_{\epsilon}D(T^{\epsilon})$, we can get a diagram in $ D(T)$ by forgetting the sign $\epsilon$, and conversely, for a diagram
        $D\in D(T)$, there is a unique way of assigning a sign $\epsilon$ satisfying the compatibility condition. So we have
        \begin{align}
            D(T)=\coprod_{\epsilon}D(T^{\epsilon}).\label{form7}
        \end{align}

    }
        \par{
        It will be helpful illustrating $D(T)$ or $D(T^{\epsilon})$ using graphs. We give some examples and non-examples of $D(T)$ and $D(T^{\epsilon})$. Let $T=(a_1,b_1)(a_2,b_2)(a_3,b_3)(a_4,b_4)$.
        Then the following is a diagram over $T$
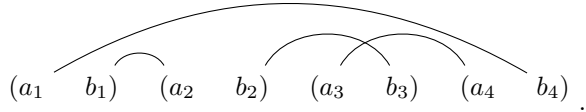
        $$
        \begin{tikzpicture}
                \node (p1) at (-3,0) {$(a_1$};
			    \node (p2) at (-2,0) {$b_1)$};
			    \node (p3) at (-1,0) {$(a_2$};
                \node (p4) at (0,0) {$b_2)$};
                \node (p5) at (1,0) {$(a_3$};
                \node (p6) at (2,0) {$b_3)$};
                \node (p7) at (3,0) {$(a_4$};
			    \node (p8) at (4,0) {$b_4)$};
                \draw (p2)to [out=60,in=120](p3);
                \draw (p1)to [out=30,in=150](p8);
                \draw (p4)to [out=60,in=120](p6);
                \draw (p5)to [out=60,in=120](p7);
        \end{tikzpicture}.
        $$
        And there is a unique way of adding the sign:
        $$
        
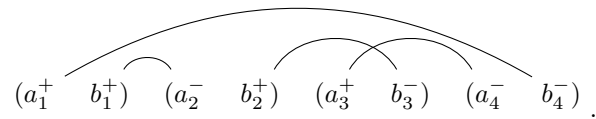
\begin{tikzpicture}
                \node (p1) at (-3,0) {$(a_1^{+}$};
			    \node (p2) at (-2,0) {$b_1^{+})$};
			    \node (p3) at (-1,0) {$(a_2^{-}$};
                \node (p4) at (0,0) {$b_2^{+})$};
                \node (p5) at (1,0) {$(a_3^{+}$};
                \node (p6) at (2,0) {$b_3^{-})$};
                \node (p7) at (3,0) {$(a_4^{-}$};
			    \node (p8) at (4,0) {$b_4^{-})$};
                \draw (p2)to [out=60,in=120](p3);
                \draw (p1)to [out=30,in=150](p8);
                \draw (p4)to [out=60,in=120](p6);
                \draw (p5)to [out=60,in=120](p7);
        \end{tikzpicture}.
        $$
        But the followings are not diagrams over $T$:
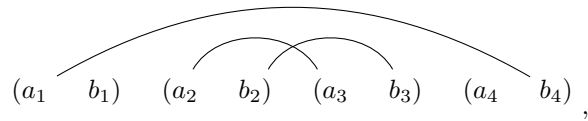
        $$
        \begin{tikzpicture}
                \node (p1) at (-3,0) {$(a_1$};
			    \node (p2) at (-2,0) {$b_1)$};
			    \node (p3) at (-1,0) {$(a_2$};
                \node (p4) at (0,0) {$b_2)$};
                \node (p5) at (1,0) {$(a_3$};
                \node (p6) at (2,0) {$b_3)$};
                \node (p7) at (3,0) {$(a_4$};
			    \node (p8) at (4,0) {$b_4)$};
                \draw  (p3) to [out=60,in=120] (p5);
                \draw  (p6) to [out=120,in=60] (p4);
                \draw  (p1) to [out=30,in=150] (p8);
        \end{tikzpicture},
        
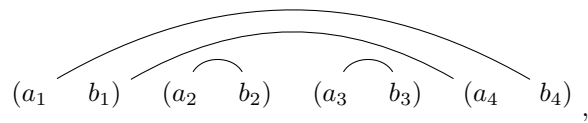
$$
        $$
        \begin{tikzpicture}
                \node (p1) at (-3,0) {$(a_1$};
			    \node (p2) at (-2,0) {$b_1)$};
			    \node (p3) at (-1,0) {$(a_2$};
                \node (p4) at (0,0) {$b_2)$};
                \node (p5) at (1,0) {$(a_3$};
                \node (p6) at (2,0) {$b_3)$};
                \node (p7) at (3,0) {$(a_4$};
			    \node (p8) at (4,0) {$b_4)$};
                \draw  (p5) to [out=60,in=120] (p6);
                \draw  (p4) to [out=120,in=60] (p3);
                \draw  (p1) to [out=30,in=150] (p8);
                \draw  (p2) to [out=30,in=150] (p7);
        \end{tikzpicture},
        
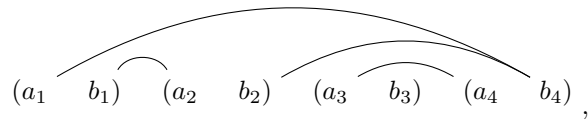
$$
        $$
        \begin{tikzpicture}
                \node (p1) at (-3,0) {$(a_1$};
			    \node (p2) at (-2,0) {$b_1)$};
			    \node (p3) at (-1,0) {$(a_2$};
                \node (p4) at (0,0) {$b_2)$};
                \node (p5) at (1,0) {$(a_3$};
                \node (p6) at (2,0) {$b_3)$};
                \node (p7) at (3,0) {$(a_4$};
			    \node (p8) at (4,0) {$b_4)$};
                \draw (p1) to [out=30,in=150] (p8);
                \draw (p4) to [out=30,in=150] (p8);
                \draw (p2) to [out=60,in=120] (p3);
                \draw (p5) to [out=30,in=150] (p7);

        \end{tikzpicture},
        $$
        because they violate the conditions (1),(2), and (3) in the definition respectively. The following
        $$
        
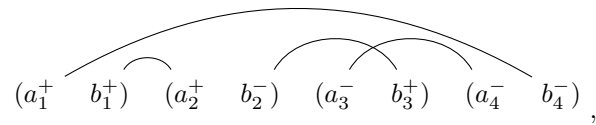
\begin{tikzpicture}
                \node (p1) at (-3,0) {$(a_1^{+}$};
			    \node (p2) at (-2,0) {$b_1^{+})$};
			    \node (p3) at (-1,0) {$(a_2^{+}$};
                \node (p4) at (0,0) {$b_2^{-})$};
                \node (p5) at (1,0) {$(a_3^{-}$};
                \node (p6) at (2,0) {$b_3^{+})$};
                \node (p7) at (3,0) {$(a_4^{-}$};
			    \node (p8) at (4,0) {$b_4^{-})$};
                \draw (p2)to [out=60,in=120](p3);
                \draw (p1)to [out=30,in=150](p8);
                \draw (p4)to [out=60,in=120](p6);
                \draw (p5)to [out=60,in=120](p7);
        \end{tikzpicture},
        $$
        is a diagram over $T$, but it is not compatible with the sign $\epsilon=(++)(+-)(-+)(--)$.
    }
    \par{
        We will establish a map from $D(T)$ to $DR(T)$ in two steps.
    }
    \par{
         First, for a diagram $D=\{V,E\}\in D(T)$, we can get a new graph by identifying two vertices $a_i,b_i$  in the same pair $(a_i,b_i)$, $\forall i=1,\cdots,n$. The defining condition of $D(T)$ implies that this graph is always disjoint unions of cycle graphs.
    }
    \par{
        For example, when $T=(a_1,b_1)\cdots(a_6,b_6)$£¬ and $D(T)$ is graphically represented by:
        
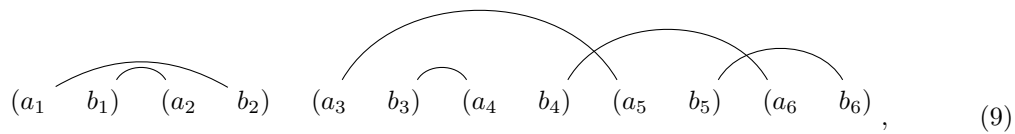
\begin{align}
        \begin{tikzpicture}
                \node (p1) at (-3,0) {$(a_1$};
			    \node (p2) at (-2,0) {$b_1)$};
			    \node (p3) at (-1,0) {$(a_2$};
                \node (p4) at (0,0) {$b_2)$};
                \node (p5) at (1,0) {$(a_3$};
                \node (p6) at (2,0) {$b_3)$};
                \node (p7) at (3,0) {$(a_4$};
			    \node (p8) at (4,0) {$b_4)$};
                \node (p9) at (5,0) {$(a_5$};
                \node (p10) at (6,0) {$b_5)$};
                \node (p11) at (7,0) {$(a_6$};
			    \node (p12) at (8,0) {$b_6)$};
                \draw (p1) to [out=30,in=150] (p4);
                \draw (p2) to [out=60,in=120] (p3);
                \draw (p5) to [out=60,in=120] (p9);
                \draw (p10) to [out=60,in=120] (p12);
                \draw (p8) to [out=60,in=120] (p11);
                \draw (p6) to [out=60,in=120] (p7);
        \end{tikzpicture},
        \label{form18}
        \end{align}
        the corresponding graph is given by
        $$
            
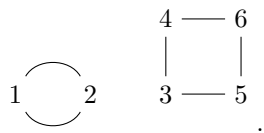
\begin{tikzpicture}
			\begin{centering}
			\node (p1) at (-1,0) {$1$};
			\node (p2) at (0,0) {$2$};
			\node (p3) at (1,0) {$3$};
            \node (p4) at (1,1) {$4$};
            \node (p5) at (2,1) {$6$};
            \node (p6) at (2,0) {$5$};
            \draw (p1) to[out=60,in=120]  (p2);
            \draw (p2) to[out=-120,in=-60] (p1);
            \draw(p3) edge(p4);
            \draw(p4) edge(p5);
            \draw(p5) edge(p6);
            \draw(p6) edge(p3);
            \end{centering}
            \end{tikzpicture}.
        $$
    }
    \par{
        Next, we add an orientation on each cycle to get a directed graph. For a cycle with vertices labelled by $i_1,\cdots,i_t$, in which $i_1$ be the smallest number, when $t\geq 3$, there are two possibilities adding the orientation. We choose the edge $\{a_{i_1},b_j\}\in E$ (or $\{a_{i_1},a_j\}\in E$), and take the orientation in the direction from $i_1$ to $j$.
    }
    \par{
        For example, for diagram (\ref{form18}),
        the corresponding directed graph is
                $$
            \begin{tikzpicture}
			\begin{centering}
			\node (p1) at (-1,0) {$1$};
			\node (p2) at (0,0) {$2$};
			\node (p3) at (1,0) {$3$};
            \node (p4) at (1,1) {$4$};
            \node (p5) at (2,1) {$6$};
            \node (p6) at (2,0) {$5$};
            \draw [->](p1) to[out=30,in=150]  (p2);
            \draw [->](p2) to[out=-150,in=-30] (p1);
            \draw(p3) edge[->](p6);
            \draw(p6) edge[->](p5);
            \draw(p5) edge[->](p4);
            \draw(p4) edge[->](p3);
            \end{centering}
            \end{tikzpicture},
        $$
        where $i_1=3,j=5$, in the second cycle. Noting that these graphs are essentially the same as what is called 'Virasoro graphs' in \cite{Hur}.
    }
    \par{
        For every such directed graph we can get a corresponding element in $DR(T)$ in the obvious way. For the above example, $(12)(3564)$ is the corresponding derangement. We use $\sigma_D$ to denote the derangement corresponding to $D$.
    }
    \par{
        We have the following combinatorial lemma:
        \begin{lemma}
            Given a sequence $T=(a_1,b_1)\cdots(a_n,b_n),n\geq 2$,
            the map defined above:
            $$
                D\mapsto \sigma_D.
            $$
            is surjective.
            Further, for $\sigma=(C_1),\cdots(C_s)\in DR(T)$, there are exactly $2^{n-s}$ diagrams $D\in D(T)$ with $\sigma_D=\sigma$.
        \end{lemma}
    }
    \par{
        We also consider some operations on the diagram. Suppose $T=(a_1,b_1)\cdots(a_n,b_n)$, $D=(V,E)\in D(T)$, $e=\{a_i,b_{i+1}\}\in E$, then we can form a new diagram $D_1=(V_1,E_1)$ by deleting $e$. More precisely the new sequence is: $$T_1=(a_1,b_1)\cdots(a_{i-1},b_{i-1})(a_{i+1},b_i)(a_{i+2},b_{i+2})\cdots(a_n,b_n),$$
        and this is described graphically by:
        $$
            \begin{tikzpicture}
                \node (p1) at (-4.5,-0.5) {$D_1$};
			    \node (p13) at (-6,0) {$\cdots$};
                \node (p14) at (-5,0) {$(a_{i+1}$};
                \node (p17) at (-4,0) {$b_{i})$};
			    \node (p18) at (-3,0) {$\cdots$};
                \draw  (p14) to [out=60,in=120] (p18);
                \draw  (p13) to [out=60,in=120] (p17);
                \node (p12) at (-2.5,0) {$=$};
                \node (p2) at (0.5,-.5) {$D$};
			    \node (p23) at (-2,0) {$\cdots$};
                \node (p24) at (-1,0) {$(a_i$};
                \node (p25) at (0,0) {$b_i)$};
                \node (p26) at (1,0) {$(a_{i+1}$};
                \node (p27) at (2,0) {$b_{i+1})$};
			    \node (p28) at (3,0) {$\cdots$};
                \draw  (p23) to [out=60,in=120] (p25);
                \draw  (p24) to [out=60,in=120] (p27);
                \draw  (p26) to [out=60,in=120] (p28);
                \node (p22) at (4,0) {$-$};
                \node (p3) at (5.5,-.5) {$e$};
                \node (p34) at (5,0) {$a_i$};
                \node (p37) at (6,0) {$b_{i+1}$};
                \draw  (p34) to [out=60,in=120] (p37);
            \end{tikzpicture}.
        $$
        We simply write $D_1=D-e$, or equivalently $D=D_1+e$.
}
\par{
        Similarly, if $e_1=\{a_i,a_{i+1}\},e_2=\{b_i,b_{i+1}\}\in E$, we can form a new diagram $D_2=(V_2,E_2)$ by deleting both $e_1,e_2$:
        $$
            T_2=(a_1,b_1)\cdots(a_{i-1},b_{i-1})(a_{i+2},b_{i+2})\cdots(a_n,b_n),
        $$
        which is described graphically as
        $$
            \begin{tikzpicture}
                \node (p1) at (-5,0) {$\cdots$};
                \node (p12) at (-4.5,-.5) {$D_2$};
                \node (p2) at (-4,0) {$\cdots$};
                \node (p23)at (-3,0) {$=$};
			    \node (p3) at (-2,0) {$\cdots$};
                \node (p4) at (-1,0) {$(a_i$};
                \node (p5) at (0,0) {$b_i)$};
                \node (p6) at (1,0) {$(a_{i+1}$};
                \node (p56) at (1.5,-0.5) {$D$};
                \node (p7) at (2,0) {$b_{i+1})$};
                \node (p78)at (3,0) {$-$};
                \node (p8) at (4,0) {$a_i$};
                \node (p89) at (4.5,-0.5) {$e_1$};
                \node (p9) at (7,0) {$b_i$};
                \node (p78)at (6,0) {$-$};
                \node (p10) at (5,0) {$a_{i+1}$};
                \node (p1011) at (7.5,-0.5) {$e_2$};
                \node (p11) at (8,0) {$b_{i+1}$};
                \draw  (p4) to [out=60,in=120] (p6);
                \draw  (p5) to [out=60,in=120] (p7);
                \draw  (p8) to [out=60,in=120] (p10);
                \draw  (p9) to [out=60,in=120] (p11);
            \end{tikzpicture}.
        $$
        We also write $D_2=D-e_1-e_2$ or $D=D_2+e_1+e_2$. Noting that $D_2=\emptyset$ may happen.
    }
}
\par{
    Later on we use capital letters $Z,W$ to denote the set of variables in the formal power series $L_{a_1,b_1}(z_1,w_1),\cdots,L_{a_n,b_n}(z_n,w_n)$:
    $$
        Z=\{z_1,\cdots,z_n\},\;W=\{w_1,\cdots,w_n\}.
    $$
    In (\ref{form5}) and (\ref{form6}) it is seen that $z_i$ is always associated to $a_i$ and $w_i$ is associated to $b_i$.
}
\par{
    We need to define some functions related to the graphical data.
    For a diagram $D=(V,E)$ over $T=(a_1,b_1)\cdots(a_n,b_n)$ and
    for an edge $e=\{a_i,b_j\}\in E$, define
    \begin{align*}
        K(e;Z,W)\stackrel{def}{=}\frac{1}{(z_i-w_j)^2},\;\;
        Q(e;Z,W)\stackrel{def}{=}\frac{\frac{1}{2}(a_i,b_j)}{(z_i-w_j)^2},
    \end{align*}
    and similarly for an edge $e=\{a_i,a_j\}$, define
    \begin{align*}
        K(e;Z,W)\stackrel{def}{=}\frac{1}{(z_i-z_j)^2},\;\;
        Q(e;Z,W)\stackrel{def}{=}\frac{\frac{1}{2}(a_i,a_j)}{(z_i-z_j)^2},
    \end{align*}
    and for $e=\{b_i,b_j\}$:
    \begin{align*}
        K(e;Z,W)\stackrel{def}{=}\frac{1}{(w_i-w_j)^2},\;\;
        Q(e;Z,W)\stackrel{def}{=}\frac{\frac{1}{2}(b_i,b_j)}{(w_i-w_j)^2},
    \end{align*}
    and we also define the functions $\Gamma(D),R(D;Z,W)$:
    \begin{align}
        \Gamma(D)\stackrel{def}{=}\prod_{\{u,v\}\in E}(u,v),\notag\;\;\;\;
        R(D;Z,W)\stackrel{def}{=}r^{c(\sigma_D)}\prod_{e\in E} Q(e;Z,W).\notag
    \end{align}
    Observing that $R(D;Z,W)$ can also be written as
    \begin{align}
        R(D;Z,W)\stackrel{def}{=}\Gamma(D)r^{c(\sigma_D)}\prod_{e\in E} K(e;Z,W),\label{form14}
    \end{align}
    through a direct computation. When $D=\emptyset$ we define $$R(\emptyset;Z,W)=1$$ by convention.
}

\section{Proof of Theorem 1}
\par{
    In this section we give the proof of Theorem 1. To compute the correlation function
    $$
        \langle 1^{'},L_{a_1,b_1}(z_1)\cdots L_{a_n,b_n}(z_n)1\rangle,
    $$
    we need to compute the following ``correlation function" of formal power series:
    $$
        \langle 1^{'},L_{a_1,b_1}(z_1,w_1)\cdots L_{a_n,b_n}(z_n,w_n)1\rangle.
    $$
By (\ref{form17}) we have
    \begin{align}
        \langle 1^{'},L_{a_1,b_1}(z_1,w_1)\cdots L_{a_n,b_n}(z_n,w_n)1\rangle=\sum_{\epsilon}\langle 1^{'},L^{\epsilon_1,\delta_1}_{a_1,b_1}(z_1,w_1)\cdots L^{\epsilon_n,\delta_n}_{a_n,b_n}(z_n,w_n)1\rangle,\notag
    \end{align}
    where $\epsilon=(\epsilon_1,\delta_1)\cdots (\epsilon_n,\delta_n)$ runs over all $4^n$ possible signs over $T$.
    We need to compute each term in the summation of the right hand side.
}
\par{
    When $D(T^{\epsilon})=\emptyset$, it is easy to see that
    $$
        \langle 1^{'},L^{\epsilon_1,\delta_1}_{a_1,b_1}(z_1,w_1)\cdots L^{\epsilon_n,\delta_n}_{a_n,b_n}(z_n,w_n)1\rangle=0.
    $$
 We have the following key lemma:
    \begin{lemma}
        For signed sequence $T^{\epsilon}=(a_1^{\epsilon_1},b_1^{\delta_1})\cdots(a_n^{\epsilon_n},b_n^{\delta_n})$, we have
        \begin{align}
            \langle 1^{'},L^{\epsilon_1,\delta_1}_{a_1,b_1}(z_1,w_1),\cdots,L^{\epsilon_n,\delta_n}_{a_n,b_n}(z_n,w_n)1\rangle=\sum_{D\in D(T^{\epsilon})} R(D;Z,W),\label{form10}
        \end{align}
         and the right hand side is $0$ if $D(T^{\epsilon})=\emptyset$ by convention.
    \end{lemma}

}

\par{
     \textbf{Proof.} Let $T^{\epsilon}(n,k)$ denote the set of signed sequence $T^{\epsilon}$ which has $n$ pairs, and there are $k$ consecutive  pairs of form $(a^{+},b^{+})$ appear in the leftmost of $T^{\epsilon}$.
    For example, for $$T^{\epsilon}=(a_1^{+},b_1^{+})(a_2^{+},b_2^{+})(a_3^{-},b_3^{+})(a_4^{+},b_4^{-})(a_5^{-},b_5^{-})(a_6^{-},b_6^{-}),$$
    then we have  $T^{\epsilon}\in T^{\epsilon}(6,2)$.
}
\par{We prove this by induction on $n+k$.
When $n=k=0$, obviously $T^{\epsilon}=\emptyset$ and both sides equal to $1$; When $k=0,n>0$, or $n<2k$ it is also easy to see that both sides equal to $0$.
}
    \par{
       So it is enough to show that the correctness for $T^{\epsilon}\in T^{\epsilon}(n,k-1)\cup T^{\epsilon}(n-1,k)\cup T^{\epsilon}(n-1,k-1)\cup T^{\epsilon}(n-2,k-1)$ will imply the correctness for $T^{\epsilon}\in T^{\epsilon}(n,k)$, and we only need to consider the case when $n\geq 1$ and $n\geq 2k$.
    }
    \par{
        For $T^{\epsilon}\in T^{\epsilon}(n,k)$, suppose $T^{\epsilon}= \cdots(a_k^{+},b^{+}_{k})\cdots$, then there are two cases:
    }
    \par{
        \begin{description}
            \item[Case 1] The signed sequence $T^{\epsilon} =\cdots(a_k^{+},b_k^{+})(a_{k+1}^{-},b_{k+1}^{+})\cdots$, or $\cdots(a_k^{+},b_k^{+})(a_{k+1}^{+},b_{k+1}^{-})\cdots$.
            \item[Case 2] The signed sequence $T^{\epsilon}=\cdots(a_k^{+},b_k^{+})(a_{k+1}^{-},b_{k+1}^{-})\cdots$.
        \end{description}
        Later on we always use ``$\cdots$" to denote the part that is the same as the original $T^{\epsilon}$.
    }
    \par{
         \textit{Case 1}: For the first case, by (\ref{form17}) we only need to consider the subcase when $T^{\epsilon}=\cdots(a_k^{+},b_k^{+})(a_{k+1}^{-},b_{k+1}^{+})\cdots$. Using commutation relation (\ref{form6}), we see that
        \begin{align}
            &\langle 1^{'}, \cdots L_{a_k,b_k}^{++}(z_k,w_k)L^{-+}_{a_{k+1},b_{k+1}}(z_{k+1},w_{k+1}) \cdots 1\rangle\notag\\
            =&\langle 1^{'}, \cdots L^{-+}_{a_{k+1},b_{k+1}}(z_{k+1},w_{k+1})L^{++}_{a_k,b_k}(z_k,w_k) \cdots 1\rangle\notag\\&+\langle 1, \cdots [L^{++}_{a_k,b_k}(z_k,w_k),L^{-+}_{a_{k+1},b_{k+1}}(z_{k+1},w_{k+1})] \cdots 1\rangle\notag\\
            =&\langle 1^{'}, \cdots L^{-+}_{a_{k+1},b_{k+1}}(z_{k+1},w_{k+1})L^{++}_{a_k,b_k}(z_k,w_k) \cdots 1\rangle\notag\\
            &+\frac{1}{2}(b_k,a_{k+1})(w_k-z_{k+1})^{-2} \langle 1^{'}, \cdots L_{a_k,b_{k+1}}^{++}(z_k,w_{k+1}) \cdots 1\rangle\notag\\&+\frac{1}{2}(a_k,a_{k+1})(z_k-z_{k+1})^{-2}\langle 1^{'} \cdots L_{b_k,b_{k+1}}^{++}(w_k,w_{k+1}) \cdots 1\rangle.\label{form19}
        \end{align}
    }
    \par{
        Let
        $$
            {T_1}^{\epsilon_1}=\cdots(a_{k+1}^{-},b_{k+1}^{+})(a_{k}^{+},b_{k}^{+})\cdots,\;\;
            T_2^{\epsilon_2}=\cdots(a_k^{+},b_{k+1}^{+})\cdots,\;\;
            T_3^{\epsilon_3}=\cdots(b_k^{+},b_{k+1}^{+})\cdots.
        $$
        Observe that $D(T^{\epsilon})$ can be decomposed into three disjoint subsets $D(T_1^{\epsilon_1})$, $D(T_2^{\epsilon_2}),D(T_3^{\epsilon_3})$:
        \begin{enumerate}
            \item The case when there is no edge connecting $(a_k^{+},b_k^{+})$ and $(a_{k+1}^{-},b_{k+1}^{+})$:
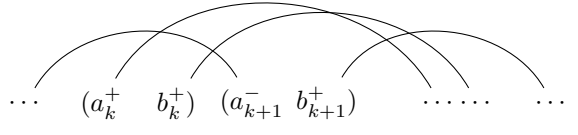
            $$
            \begin{tikzpicture}

			    \node (p2) at (-2,0) {$\cdots$};
			    \node (p3) at (-1,0) {$(a_k^{+}$};
                \node (p4) at (0,0) {$b_k^{+})$};
                \node (p5) at (1,0) {$(a_{k+1}^{-}$};
                \node (p6) at (2,0) {$b_{k+1}^{+})$};
                \node (p7) at (3.5,0) {$\cdots$};
			    \node (p8) at (4,0) {$\cdots$};
                \node (p9) at (5,0) {$\cdots$};
                \draw  (p6) to [out=60,in=120] (p9);
                \draw  (p2) to [out=60,in=120] (p5);
                \draw  (p3) to [out=60,in=120] (p7);
                \draw  (p4) to [out=60,in=120] (p8);
            \end{tikzpicture}.
        $$
            \item The case when there is one edge connecting $(a_k^{+},b_k^{+})$ and $(a_{k+1}^{-},b_{k+1}^{+})$, and there are two subcases:
            $$
            
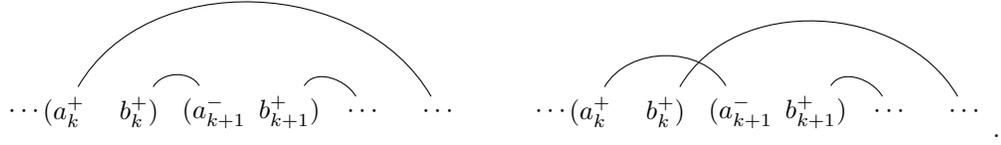
\begin{tikzpicture}
                \node (p20) at (0.5,0) {$\cdots$};
                \node (p14) at (1,0) {$(a_k^{+}$};
                \node (p15) at (2,0) {$b_k^{+})$};
                \node (p16) at (3,0) {$(a_{k+1}^{-}$};
                \node (p17) at (4,0) {$b_{k+1}^{+})$};
			    \node (p18) at (5,0) {$\cdots$};
                \node (p19) at (6,0) {$\cdots$};
                \draw  (p15) to [out=60,in=120] (p19);
                \draw  (p14) to [out=60,in=120] (p16);
                \draw  (p17) to [out=60,in=120] (p18);
                \node (p21) at (-6.5,0) {$\cdots$};
                \node (p4) at (-6,0) {$(a_k^{+}$};
                \node (p5) at (-5,0) {$b_k^{+})$};
                \node (p6) at (-4,0) {$(a_{k+1}^{-}$};
                \node (p7) at (-3,0) {$b_{k+1}^{+})$};
			    \node (p8) at (-2,0) {$\cdots$};
                \node (p9) at (-1,0) {$\cdots$};
                \draw  (p5) to [out=60,in=120] (p6);
                \draw  (p4) to [out=60,in=120] (p9);
                \draw  (p7) to [out=60,in=120] (p8);
            \end{tikzpicture}.
        $$
        \end{enumerate}
        By our notation of diagram operations defined in section 3, we have:
        \begin{align}
            D(T^{\epsilon})=(D(T_1^{\epsilon_1}))\coprod(D(T_2^{\epsilon_2})+ e_1)\coprod(D(T_3^{\epsilon_3})+e_2),\label{form11}
        \end{align}
        where $e_1=\{b^{+}_k,a^{-}_{k+1}\},e_2=\{a^{+}_k,a^{-}_{k+1}\}$.
    }
    \par{
        Noting that $T_1^{\epsilon_1}\in T^{\epsilon}(n,k-1), T_2^{\epsilon_2},T_3^{\epsilon_3}\in T^{\epsilon}(n-1,k)$, we rewrite (\ref{form19}), and then by the induction hypothesis we have:
        \begin{align*}
            &\langle 1^{'},L^{\epsilon_1,\delta_1}_{a_1,b_1}(z_1,w_1),\cdots,L^{\epsilon_n,\delta_n}_{a_n,b_n}(z_n,w_n)1\rangle\\
            =&\sum_{D\in D(T_1^{\epsilon_1})} R(D;Z,W)\\
            &+Q(e_1;Z,W)\sum_{D\in D(T_2^{\epsilon_2})}R(D;Z,W)
            +Q(e_2;Z,W)\sum_{D\in D(T_3^{\epsilon_3})}R(D;Z,W)\\
            =&\sum_{D\in D(T_1^{\epsilon_1})} R(D;Z,W)\\
            &+\sum_{D\in D(T_2^{\epsilon_2})+e_1}R(D;Z,W)
            +\sum_{D\in D(T_3^{\epsilon_3})+e_2}R(D;Z,W)
            =\sum_{D\in D(T^{\epsilon})} R(D;Z,W).
        \end{align*}
        So the proposition holds for this case.
    The third step is by observing that if $D_1=D-e$, then $c(\sigma_D)=c(\sigma_{D_1})$, so we have
        \begin{align}
            R(D;Z,W)=Q(e;Z,W)R(D_1;Z,W).\label{form8}
        \end{align}
        and the last step is given by (\ref{form11}).
}

    \par{
        \textit{Case 2}: The arguments here are similar but there are more subcases. Formula (\ref{form5}) tells that
        \begin{align}
            &\langle 1, \cdots L_{a_k,b_k}^{++}(z_k,w_k)L^{--}_{a_{k+1},b_{k+1}}(z_{k+1},w_{k+1}) \cdots 1\rangle\notag\\
            =&\langle 1, \cdots L^{--}_{a_{k+1},b_{k+1}}(z_{k+1},w_{k+1})L^{++}_{a_k,b_k}(z_k,w_k) \cdots 1\rangle\notag\\&+\langle 1, \cdots [L^{++}_{a_k,b_k}(z_k,w_k),L^{--}_{a_{k+1},b_{k+1}}(z_{k+1},w_{k+1})] \cdots 1\rangle\notag\\
            =&\langle 1, \cdots L^{--}_{a_{k+1},b_{k+1}}(z_{k+1},w_{k+1})L^{++}_{a_k,b_k}(z_k,w_k) \cdots 1\rangle\notag\\
            &+\frac{1}{2}(a_k,b_{k+1})(z_k-w_{k+1})^{-2}\langle 1, \cdots L_{a_{k+1},b_k}^{-+}(z_{k+1},w_k) \cdots 1\rangle\notag\\&+\frac{1}{2}(a_k,a_{k+1})(z_k-z_{k+1})^{-2}\langle 1, \cdots L_{b_{k+1},b_k}^{-+}(w_{k+1},w_k) \cdots 1\rangle\notag\\
            &+\frac{1}{2}(b_k,b_{k+1})(w_k-w_{k+1})^{-2} \langle1, \cdots L_{a_{k+1},a_k}^{-+}(z_{k+1},z_k) \cdots 1\rangle\notag\\&+\frac{1}{2}(b_k,a_{k+1})(w_k-z_{k+1})^{-2} \langle1, \cdots L_{b_{k+1},a_k}^{-+}(w_{k+1},z_k) \cdots 1\rangle\notag\\
            &+\frac{r}{4}(a_k,a_{k+1})(z_k-z_{k+1})^{-2}(b_k,b_{k+1})(w_k-w_{k+1})^{-2}\langle 1, \cdots  \cdots 1\rangle\notag\\
            &+\frac{r}{4}(a_k,b_{k+1})(z_k-w_{k+1})^{-2}(a_{k+1},b_k)(w_k-z_{k+1})^{-2}\langle 1, \cdots  \cdots 1\rangle.\label{form20}
        \end{align}
    }
\par{
    $D(T^{\epsilon})$ can also be decomposed into three disjoint subsets:
    \begin{enumerate}
        \item The case when there is no edge connecting $(a_k^{+},b_k^{+})$ and  $(a_{k+1}^{-},b_{k+1}^{-})$, then we can exchange these two pairs:
        $$
            
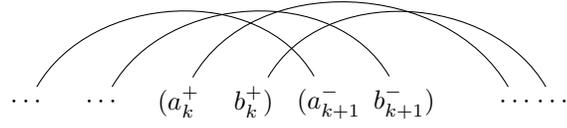
\begin{tikzpicture}
                \node (p1) at (-3,0) {$\cdots$};
			    \node (p2) at (-2,0) {$\cdots$};
			    \node (p3) at (-1,0) {$(a_k^{+}$};
                \node (p4) at (0,0) {$b_k^{+})$};
                \node (p5) at (1,0) {$(a_{k+1}^{-}$};
                \node (p6) at (2,0) {$b_{k+1}^{-})$};
                \node (p7) at (3.5,0) {$\cdots$};
			    \node (p8) at (4,0) {$\cdots$};
                \draw  (p1) to [out=60,in=120] (p5);
                \draw  (p2) to [out=60,in=120] (p6);
                \draw  (p3) to [out=60,in=120] (p7);
                \draw  (p4) to [out=60,in=120] (p8);
            \end{tikzpicture}.
        $$
        \item The case when there is exactly one edge connecting $(a_k^{+},b_k^{+})$ and  $(a_{k+1}^{-},b_{k+1}^{-})$, and there are four subcases:
        $$
            
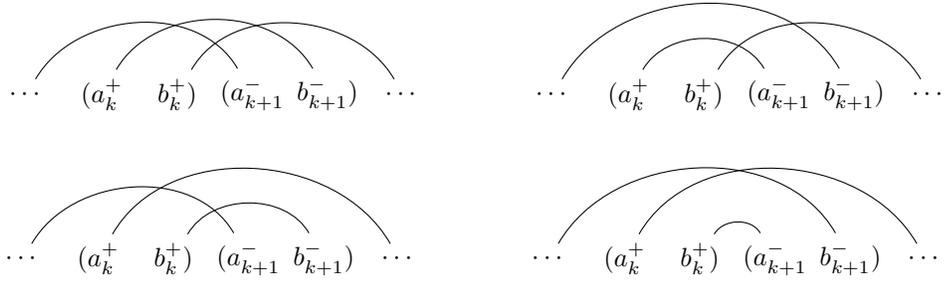
\begin{tikzpicture}
			    \node (p13) at (-6,0) {$\cdots$};
                \node (p14) at (-5,0) {$(a_k^{+}$};
                \node (p15) at (-4,0) {$b_k^{+})$};
                \node (p16) at (-3,0) {$(a_{k+1}^{-}$};
                \node (p17) at (-2,0) {$b_{k+1}^{-})$};
			    \node (p18) at (-1,0) {$\cdots$};
                \draw  (p13) to [out=60,in=120] (p16);
                \draw  (p14) to [out=60,in=120] (p17);
                \draw  (p15) to [out=60,in=120] (p18);
			    \node (p3) at (1,0) {$\cdots$};
                \node (p4) at (2,0) {$(a_k^{+}$};
                \node (p5) at (3,0) {$b_k^{+})$};
                \node (p6) at (4,0) {$(a_{k+1}^{-}$};
                \node (p7) at (5,0) {$b_{k+1}^{-})$};
			    \node (p8) at (6,0) {$\cdots$};
                \draw  (p3) to [out=60,in=120] (p7);
                \draw  (p4) to [out=60,in=120] (p6);
                \draw  (p5) to [out=60,in=120] (p8);
            \end{tikzpicture}
        $$
        $$
            \begin{tikzpicture}
			    \node (p13) at (-6,0) {$\cdots$};
                \node (p14) at (-5,0) {$(a_k^{+}$};
                \node (p15) at (-4,0) {$b_k^{+})$};
                \node (p16) at (-3,0) {$(a_{k+1}^{-}$};
                \node (p17) at (-2,0) {$b_{k+1}^{-})$};
			    \node (p18) at (-1,0) {$\cdots$};
                \draw  (p15) to [out=60,in=120] (p17);
                \draw  (p14) to [out=60,in=120] (p18);
                \draw  (p13) to [out=60,in=120] (p16);
			    \node (p3) at (1,0) {$\cdots$};
                \node (p4) at (2,0) {$(a_k^{+}$};
                \node (p5) at (3,0) {$b_k^{+})$};
                \node (p6) at (4,0) {$(a_{k+1}^{-}$};
                \node (p7) at (5,0) {$b_{k+1}^{-})$};
			    \node (p8) at (6,0) {$\cdots$};
                \draw  (p5) to [out=60,in=120] (p6);
                \draw  (p4) to [out=60,in=120] (p8);
                \draw  (p3) to [out=60,in=120] (p7);
            \end{tikzpicture}.
        $$
        \item The case when there are exactly two edges connecting $(a_k^{+}£¬b_k^{+})$ and $(a_{k+1}^{-}£¬b_{k+1}^{-})$, and there are two subcases£º
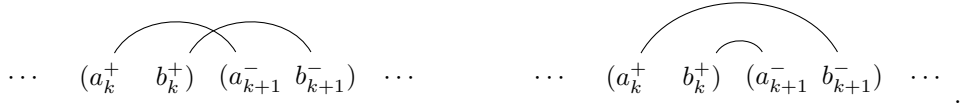
        $$
            \begin{tikzpicture}
			    \node (p13) at (-6,0) {$\cdots$};
                \node (p14) at (-5,0) {$(a_k^{+}$};
                \node (p15) at (-4,0) {$b_k^{+})$};
                \node (p16) at (-3,0) {$(a_{k+1}^{-}$};
                \node (p17) at (-2,0) {$b_{k+1}^{-})$};
			    \node (p18) at (-1,0) {$\cdots$};
                \draw  (p14) to [out=60,in=120] (p16);
                \draw  (p15) to [out=60,in=120] (p17);
			    \node (p3) at (1,0) {$\cdots$};
                \node (p4) at (2,0) {$(a_k^{+}$};
                \node (p5) at (3,0) {$b_k^{+})$};
                \node (p6) at (4,0) {$(a_{k+1}^{-}$};
                \node (p7) at (5,0) {$b_{k+1}^{-})$};
			    \node (p8) at (6,0) {$\cdots$};
                \draw  (p4) to [out=60,in=120] (p7);
                \draw  (p5) to [out=60,in=120] (p6);
            \end{tikzpicture}.
        $$
    \end{enumerate}
}
\par{
    Let
    \begin{align*}
        &T_1^{\epsilon_1}=\cdots(a_{k+1}^{-},b_{k+1}^{-})(a_k^{+},b_k^{+})\cdots,\;T_3^{\epsilon_3}=\cdots \cdots\\
        &T_{21}^{\epsilon_{21}}=\cdots(a_{k+1}^{-},b_k^{+})\cdots,\;
        T_{22}^{\epsilon_{22}}=\cdots(b_{k+1}^{-},b_k^{+})\cdots,\\
        &T_{23}^{\epsilon_{23}}=\cdots(a_{k+1}^{-},a_k^{+})\cdots,\;
        T_{24}^{\epsilon_{24}}=\cdots(b_{k+1}^{-},a_k^{+})\cdots.
    \end{align*}
    and
    \begin{align*}
        e_{1}=\{a_k,b_{k+1}\},\;
        e_{2}=\{a_k,a_{k+1}\},\;
        e_{3}=\{b_k,b_{k+1}\},\;
        e_{4}=\{b_k,a_{k+1}\}.
    \end{align*}
    Then we also have set decomposition:
    \begin{align}
        D(T^{\epsilon})=&(D(T_1^{\epsilon_1}))\coprod(D(T_{21}^{\epsilon_{21}})+e_1)\coprod(D(T_{22}^{\epsilon_{22}})+e_2)\coprod(D(T_{23}^{\epsilon_{23}})+e_3)\coprod(D(T_{24}^{\epsilon_{24}})+e_4)\coprod\notag\\
        &(D(T_3^{\epsilon_3})+e_2+e_3)\coprod(D(T_3^{\epsilon_3})+e_1+e_4).\label{form12}
    \end{align}
    Because $T_1^{\epsilon_1}\in T^{\epsilon}(n,k-1),T_{2i}^{\epsilon_{2i}}\in T^{\epsilon}(n-1,k-1),T_3^{\epsilon_3}\in T^{\epsilon}(n-2,k-1)$, we rewrite \ref{form20}, by definition of $Q(e;Z,W)$, $R(D;Z,W)$, and then by induction hypothesis we have:
    \begin{align*}
        \langle &1^{'},L^{\epsilon_1,\delta_1}_{a_1,b_1}(z_1,w_1),\cdots,L^{\epsilon_n,\delta_n}_{a_n,b_n}(z_n,w_n)1\rangle\notag\\
        =&\sum_{D\in D(T_1^{\epsilon_1})}R(D;Z,W)\\
        &+Q(e_1;Z,W)\sum_{D\in D(T_{21}^{\epsilon_{21}})}R(D;Z,W)
        +Q(e_2;Z,W)\sum_{D\in D(T_{22}^{\epsilon_{22}})}R(D;Z,W)\\
        &+Q(e_3;Z,W)\sum_{D\in D(T_{23}^{\epsilon_{23}})}R(D;Z,W)
        +Q(e_4;Z,W)\sum_{D\in D(T_{24}^{\epsilon_{24}})}R(D;Z,W)\\
        &+(rQ(e_2;Z,W)Q(e_3;Z,W)+rQ(e_1;Z,W)Q(e_4;Z,W))\sum_{D\in D(T_3^{\epsilon_3})}R(D;Z,W)\\
        =&\sum_{D\in D(T_1^{\epsilon_1})}R(D;Z,W)+\sum_{D\in D(T_{21}^{\epsilon_{21}})+e_1}R(D;Z,W)
        +\sum_{D\in D(T_{22}^{\epsilon_{22}})+e_2}R(D;Z,W)\\
        &+\sum_{D\in D(T_{23}^{\epsilon_{23}})+e_3}R(D;Z,W)
        +\sum_{D\in D(T_{24}^{\epsilon_{24}})+e_4}R(D;Z,W)\\
        &+\sum_{D\in D(T_3^{\epsilon_3})+e_2+e_3}R(D;Z,W)+\sum_{D\in D(T_3^{\epsilon_3})+e_1+e_4}R(D;Z,W)=\sum_{D\in D(T^{\epsilon})}R(D;Z,W),
    \end{align*}
    where we also use (\ref{form8}),(\ref{form12}), and noting that for $D_2=D-e_1-e_2$, then
    $c(D_2)=c(D)-1$, and we have:
        \begin{align}
            R(D;Z,W)=rQ(e_1;Z,W)Q(e_2;Z,W)R(D_2;Z,W).
        \end{align} So we conclude the proof of (\ref{form10}). $\qed$
}
\par{
    Using (\ref{form7}) and (\ref{form10}) we see that
    \begin{align}
        \langle 1^{'},L_{a_1,b_1}(z_1,w_1)\cdots L_{a_n,b_n}(z_n,w_n)1\rangle=&\sum_{\epsilon}\langle 1^{'},L^{\epsilon_1,\delta_1}_{a_1,b_1}(z_1,w_1)\cdots L^{\epsilon_n,\delta_n}_{a_n,b_n}(z_n,w_n)1\rangle\notag\\
        =&\sum_{\epsilon}\sum_{D\in D(T^{\epsilon})}R(D;Z,W)=\sum_{D\in D(T)}R(D;Z,W).\notag
    \end{align}
    So we have:
    \begin{proposition}
    \begin{align}
        \langle 1^{'},L_{a_1,b_1}(z_1,w_1),\cdots,L_{a_n,b_n}(z_n,w_n)1\rangle=&\sum_{D\in D(T)}R(D;Z,W).\label{form13}
    \end{align}
    \end{proposition}
}
\par{
    Theorem 1 is a corollary of Proposition 2 by letting $z_1=w_1,\cdots,z_n=w_n$. Noting that for $D\in DR(T)$ we have
    \begin{align}
        \prod_{e\in E_D}K(e,Z,Z)=\Gamma(\sigma_D;Z).\label{form15}
    \end{align}
    Also note that
     $\forall \sigma\in DR(T)$, we have
        \begin{align}
            \Gamma(\sigma;T)=\sum_{D\in D(T), s.t. \sigma_D=\sigma}\Gamma(D),\label{form16}
        \end{align}
     by a direct computation.
    }
    \par{
    Continue with (\ref{form13}) we have
    \begin{align*}
        \langle 1^{'},L_{a_1,b_1}(z_1)\cdots L_{a_n,b_n}(z_n))\cdot 1\rangle
        =&\sum_{\epsilon}\langle 1^{'},L^{\epsilon_1,\delta_1}_{a_1,b_1}(z_1,w_1),\cdots,L^{\epsilon_n,\delta_n}_{a_n,b_n}(z_n,w_n)1\rangle\\
        =&\sum_{D\in D(T)}R(D;Z,Z)\label{form13}\\
        =&\sum_{D\in D(T)}(\Gamma(D)r^{c(\sigma_D)}\prod_{e\in E_D} K(e;Z,Z))\\
        =&\sum_{D\in D(T)}(\Gamma(D)\Gamma(\sigma_D;Z)r^{c(\sigma_D)})\\
        =&\sum_{\sigma\in DR(T)}(\sum_{D\in D(T),\sigma_D=\sigma}\Gamma(D)) \Gamma(\sigma;Z)r^{c(\sigma)})\\
        =&\sum_{\sigma\in DR(T)}\Gamma(\sigma;T)\Gamma(\sigma;Z)r^{c(\sigma)}.
    \end{align*}
    where we use (\ref{form14})(\ref{form15}) and (\ref{form16}). So we finally get (\ref{form1}) and the Theorem 1 is proved. $\qed$
}
\bibliographystyle{alpha}

\par{
\textsc{Department of Mathematics, The Hong Kong University of Science and Technology, Clear Water Bay, Kowloon}\\
}
\par{
\textit{Email Address}: \textbf{hzhaoab@ust.hk}
}
\end{document}